\documentclass[12pt,twoside]{article}

\usepackage{amsfonts}
\usepackage{amssymb}
\usepackage{amsmath}
\usepackage{graphicx}

\newcommand{\ignore}[1]{}

\def\@begintheorem#1#2{\par\bgroup{\sc #1\ #2. }\it\ignorespaces}
\def\@opargbegintheorem#1#2#3{\par\bgroup{\sc #1\ #2\ (#3). }\it\ignorespaces}
\def\@endtheorem{\egroup}
\newtheorem{theorem}{Theorem}[section]
\newtheorem{corollary}[theorem]{Corollary}
\newtheorem{lemma}[theorem]{Lemma}
\newtheorem{proposition}[theorem]{Proposition}
\newtheorem{example}[theorem]{Example}
\newtheorem{algorithm}[theorem]{Algorithm}
\newtheorem{definition}[theorem]{Definition}
\newcommand{\bt}[1]{\begin{theorem}\label{#1}}
\newcommand{\bc}[1]{\begin{corollary}\label{#1}}
\newcommand{\bl}[1]{\begin{lemma}\label{#1}}
\newcommand{\bp}[1]{\begin{proposition}\label{#1}}
\newcommand{\be}[1]{\begin{example}\label{#1}}
\newcommand{\ba}[1]{\begin{algorithm}\rm\label{#1}}
\newcommand{\bd}[1]{\begin{definition}\rm\label{#1}}
\newcommand{\bpr}{\noindent {\em Proof. }}
\newcommand{\et}{\end{theorem}}
\newcommand{\ec}{\end{corollary}}
\newcommand{\el}{\end{lemma}}
\newcommand{\ep}{\end{proposition}}
\newcommand{\ee}{\end{example}}
\newcommand{\ea}{\end{algorithm}}
\newcommand{\ed}{\end{definition}}
\newcommand{\epr}{{\ \vbox{\hrule\hbox{%
\vrule height1.3ex\hskip0.8ex\vrule}\hrule}}\\\par}

\def\R{\mathbb{R}}
\def\Z{\mathbb{Z}}

\def \K {{\cal K}}
\def \M {{\cal M}}
\def \G {{\cal G}}
\def \l {\langle}
\def \r {\rangle}

\def\Zinf{\Z_{\infty}}

\def \L {{\cal L}}
\def \G {{\cal G}}

\def \C {{\cal C}}
\def \Q {{\cal Q}}
\def \S {{\cal S}}
\def \P {{\cal P}}
\def \t {^{^\intercal}}
\def \l {\langle}
\def \r {\rangle}

\def \D {{\cal D}}

\def \cone {{\rm cone}}
\def \supp {{\rm supp}}

\def \diag {{\rm diag}}
\def \Diag {{\rm Diag}}

\begin{document}

\title{\bf The Quadratic Graver Cone,
Quadratic Integer Minimization, and Extensions}

\author{Jon Lee, Shmuel Onn, Lyubov Romanchuk, Robert Weismantel}

\date{}

\maketitle

\begin{abstract}
We consider the nonlinear integer programming problem of minimizing
a quadratic function over the integer points in variable
dimension satisfying a system of linear inequalities.
We show that when the Graver basis of the matrix defining the
system is given, and the quadratic function lies
in a suitable {\em dual Graver cone}, the problem can be solved in
polynomial time. We discuss the relation between this cone and the
cone of positive semidefinite matrices, and show that none contains
the other. So we can minimize in polynomial time some non-convex
and some (including all separable) convex quadrics.

We conclude by extending our results to efficient integer
minimization of multivariate polynomial functions of
arbitrary degree lying in suitable cones.
\end{abstract}

\section{Introduction}

Consider the general nonlinear integer minimization
problem in standard form,
\begin{equation}\label{NIM}
\min\,\left\{f(x)\ :\ x\in\Z^n\,,\ Ax=b\,,\ l\leq x\leq u\right\},
\end{equation}
with $A\in\Z^{m\times n}$, $b\in\Z^m$,
$l,u\in\Z_{\infty}^n$ with $\Z_{\infty}:=\Z\uplus\{\pm\infty\}$,
and $f:\R^n\rightarrow\R$.

It is well known to be NP-hard already for linear functions.
However, recently it was shown that, if the {\em Graver basis}
$\G(A)$ of $A$ is given as part of the input, then
the problem can be solved in polynomial time for the following
classes of functions. First, in \cite{DHORW}, for composite
concave functions $f(x)=g(Wx)$, with $W\in\Z^{d\times n}$,
$g:\R^d\rightarrow\R$ concave, and $d$ fixed. Second,
in \cite{HOW}, for separable convex functions
$f(x)=\sum_i f_i(x_i)$ with each $f_i$ univariate convex,
and in particular for linear functions $f(x)=w\t x$.
While the Graver basis is a complex object, it can be computed
in polynomial time from $A$ for many natural and useful
classes of matrices as demonstrated in \cite{DHORW,HOW}.
Moreover, the results of \cite{DO} imply that there is a
parameterized scheme that enables to construct increasingly
better approximations of the Graver basis of any matrix $A$ and
obtain increasingly better approximations to problem (\ref{NIM}),
see \cite{Onn} for details.

In this article we continue this line of investigation and consider
problem (\ref{NIM}) for quadratic functions
$f(x)=x\t Vx+w\t x+a$ with $V\in\R^{n\times n}$, $w\in\R^n$, and
$a\in\R$. We also discuss extensions to multivariate
polynomial functions of arbitrary degree.

We begin by noting that problem (\ref{NIM}) remains NP-hard
even if the Graver basis is part of the input and even if
the objective function is quadratic convex of rank $1$.
\bp{QuadraticHardness}
It is NP-hard to determine the optimal value of the problem
\begin{equation}\label{NPHard}
\min\left\{x\t Vx+w\t x+a\ :\ x\in\Z^n\,,
\ Ax=b\,,\ l\leq x\leq u\right\}
\end{equation}
even when $\G(A)$ is given and the function
is convex quadratic with matrix $V=vv\t$.
\ep
\bpr
Let $v\in\Z_+^n$ and $v_0\in\Z_+$ be input to the {\em subset sum}
problem of deciding if there exists $x\in\{0,1\}^n$ with
$v\t x=v_0$. Let $A:=0$ be the zero $1\times n$ matrix,
whose Graver basis $\G(A)=\{\pm{\bf 1}_i\,:\,i=1,\dots,n\}$
consists of the $n$ unit vectors and their negations.
Let $l:=0$ and $u:={\bf 1}$ be the zero and all-ones vectors
in $\Z^n$, and let $b:=0$ in $\Z^m$. Let $V:=vv\t$,
$w:=-2v_0v$, and $a:=v_0^2$. Then problem (\ref{NPHard}) becomes
\begin{equation*}
\min\left\{\left(v\t x-v_0\right)^2\ :\ x\in\{0,1\}^n\right\},
\end{equation*}
whose optimal value is $0$ if and only if
there is a subset sum, proving the claim.
\epr
This shows that to solve problem (\ref{NPHard}) in polynomial time,
even when the Graver basis is given, some restrictions on the
class of quadratic functions must be enforced.

In Section 2 we introduce the {\em quadratic Graver cone}
$\Q(A)$, which is a cone of $n\times n$
matrices defined via the Graver basis of $A$, and the
{\em diagonal Graver cone} $\D(A)$ which is the diagonal
projection of $\Q(A)$ into $\R_+^n$. We discuss some
elementary properties of these cones and their duals
$\Q^*(A)$ and $\D^*(A)$ and give some examples.

In Section 3 we prove the following algorithmic result
about the solvability of problem (\ref{NIM}) for every quadratic
function (possibly indefinite, neither convex nor concave)
whose defining matrix lies in the dual quadratic Graver cone.

\bt{QuadraticIntegerMinimization}
There is an algorithm that, given $\G(A)$,
solves the quadratic problem
\begin{equation}\label{QIM}
\min\left\{x\t Vx+w\t x+a\ :\ x\in\Z^n\,,
\ Ax=b\,,\ l\leq x\leq u\right\}
\end{equation}
in polynomial time for every integer matrix $V$ lying
in the cone $\Q^*(A)$ dual to $\Q(A)$.
\et
We point out that, in practice, the algorithm that underlies
Theorem \ref{QuadraticIntegerMinimization} can be applied
to any quadratic function. The algorithm will always stop
and output a feasible solution if one exists, which can be
used as an approximation of the optimal one. And, whenever
$V$ lies in $\Q^*(A)$, the solution produced will be true optimal.

As a special case we obtain the following
result on separable quadratic functions.
\bt{DiagonalIntegerMinimization}
There is an algorithm that, given $\G(A)$,
solves the separable problem
\begin{equation}\label{DiagonalEquation}
\min\{\sum_{i=1}^n\left(v_i x_i^2+w_ix_i+a_i\right)\ :\
x\in\Z^n\,,\ Ax=b\,,\ l\leq x\leq u\}
\end{equation}
in polynomial time for every integer vector $v$ lying in the cone
$\D^*(A)$ dual to $\D(A)$. In particular, this applies to any
convex separable quadratic, that is, with $v\in\Z_+^n$.
\et
In particular, Theorem \ref{DiagonalIntegerMinimization} enables us to
solve the problem with any linear objective function $f(x)=w\t x$,
which is the special case with $v=0$, which is always in $\D^*(A)$.

In Section 4 we proceed with a discussion of the relation between
the dual quadratic Graver cone $\Q^*(A)$ and the cone $\S_+^n$
of symmetric positive semidefinite matrices, and establish
Theorem \ref{Characterization} which provides a characterization,
in terms of their matroids only, of those matrices $A$ for which the
dual diagonal Graver cone $\D^*(A)$ strictly contains $\R_+^n$
and for which Theorem \ref{DiagonalIntegerMinimization} assures
efficient solution of problem (\ref{DiagonalEquation}) for all
separable convex as well as some nonconvex quadratic functions.

In the final Section 5 we extend our results to multivariate
polynomial functions of arbitrary degree. We define a hierarchy
of higher degree analogues $\P_k(A)$ of the quadratic Graver cone,
and show that the iterative algorithm of
Theorem \ref{QuadraticIntegerMinimization} solves the polynomial
integer minimization problem (\ref{NIM}) in polynomial time
for every degree $d$ form $f$ that lies in a cone $\K_d(A)$
defined in terms of the dual Graver cones $\P_k^*(A)$.

\bt{PolynomialIntegerMinimization}
For every fixed $d$ there is an algorithm that, given $\G(A)$, solves
\begin{equation}\label{PIM}
\min\left\{f(x)\ :\ x\in\Z^n\,,\ Ax=b\,,\ x\geq 0\right\}
\end{equation}
in polynomial time for every degree $d$ integer homogenous
polynomial $f$ in $\K_d(A)$.
\et

\section{The quadratic and diagonal Graver cones}

We begin with some notation. The inner product of two $m\times n$
matrices $U,V$ is $U\cdot V:=\sum_{i,j}U_{i,j}V_{i,j}$.
The diagonal of $n\times n$ matrix $V$ is the vector
$v:=\diag(V)\in\R^n$ defined by $v_i:=V_{i,i}$ for all $i$.
For $u\in\R^n$ we denote by $U:=\Diag(u)$ the $n\times n$
diagonal matrix with $\diag(U)=u$.
The {\em pointwise product} of vectors $g,h\in\R^n$ is the vector
$g\circ h$ in $\R^n$ with $(g\circ h)_i:=g_ih_i$ for all $i$.
Note that $g,h$ lie in the same orthant of $\R^n$ if and only
if $g\circ h\geq 0$. The {\em tensor product} of $g,h\in\R^n$
is the $n\times n$ matrix $g\otimes h=gh\t$ with
$(g\otimes h)_{i,j}:=(gh\t)_{i,j}=g_ih_j$ for all $i,j$.
We will use the notation $g\otimes h$ and $gh\t$ interchangeably
as we find appropriate. Note that for all $g,h\in\R^n$ and
$V\in\R^{n\times n}$, we have $g\circ h=\diag(g\otimes h)$
and $(g\otimes h)\cdot V=g\t Vh$.

Any quadratic function $f(x)=x\t Vx+w\t x+a$ has an equivalent
description $f(x)=x\t Ux+w\t x+a$ with $U:={1\over2}(V+V\t)$
symmetric matrix. We therefore can and will be working with
symmetric matrices which are much better behaved than arbitrary
square matrices. We denote by $\S^n\subset\R^{n\times n}$
the linear subspace of symmetric $n\times n$ matrices.
A {\em cone} is a subset $\P$ of real vector space such that
$\alpha x+\beta y\in\P$ for all $x,y\in\P$ and
$\alpha,\beta\in\R_+$. The {\em cone generated} by a set $\cal V$
of vectors is the set $\cone({\cal V})$ of nonnegative linear
combinations of finitely many vectors from $\cal V$.
In particular, $\cone(\emptyset):=\{0\}$. We will be using cones
$\D\subseteq\R^n$ of vectors and cones $\Q\subseteq\S^n$
of $n\times n$ symmetric matrices. The {\em dual} of a cone
$\D\subseteq\R^n$ and the (symmetric) {\em dual} of a cone
$\Q\subseteq\S^n$ are, respectively, the cones
$$\D^*\ :=\ \{v\in\R^n\ :\ u\t v\geq 0\,,\ u\in\D\},\quad
\Q^*\ :=\ \{V\in\S^n\ :\ U\cdot V\geq 0\,,\ U\in\Q\}.$$
Duality reverses inclusions, that is, if $\P\subseteq\K$
are cones in $\R^n$ or $\S^n$ then $\K^*\subseteq\P^*$.

We proceed with the definition of the Graver basis of an
integer matrix. The lattice of an integer $m\times n$
matrix $A$ is the set $\L(A):=\{x\in\Z^n\,:\, Ax=0\}$.
We denote by $\L^*(A)$ the set of nonzero elements in $\L(A)$.
We use a partial order $\sqsubseteq$
on $\R^n$ which extends the coordinate-wise partial order
$\leq$ on the nonnegative orthant $\R^n_+$ and is defined
as follows. For $x,y\in\R^n$ we write $x\sqsubseteq y$ and
say that $x$ is {\em conformal} to $y$ if $x\circ y\geq 0$
(that is, $x,y$ lie in the same orthant) and $|x_i|\leq |y_i|$
for all $i$. We write $x\sqsubset y$ if $x\sqsubseteq y$
and $x\neq y$. A simple extension of the classical Gordan
Lemma implies that every subset of $\Z^n$ has finitely
many $\sqsubseteq$-minimal elements.
\bd{GraverBasisDefinition}
The {\em Graver basis}\index{Graver basis} of an
integer matrix $A$ is defined to be the finite set
$\G(A)\subset\Z^n$ of $\sqsubseteq$-minimal elements
in $\L^*(A)=\{x\in\Z^n\,:\, Ax=0,\ x\neq 0\}$.
\ed

In this article we introduce the following
objects defined via the Graver basis.
\bd{QuadraticGraverCone}
The {\em quadratic Graver cone} of an integer $m\times n$ matrix $A$
is defined to be the cone $\Q(A)\subseteq\S^n$ of $n\times n$ matrices
generated by the matrices $g\otimes h+h\otimes g$ over all pairs of
distinct elements $g,h\in\G(A)$ that lie in the same orthant, that is,
$$\Q(A)\ :=\
\cone\left\{g\otimes h+h\otimes g\ :\ g,h\in\G(A)\,,\
g\neq h\,,\ g\circ h\geq 0\right\}\ \subseteq\ \S^n\ .$$
The {\em dual quadratic Graver cone}
is its (symmetric) dual $\Q^*(A)$ in $\S^n$ given by
\begin{eqnarray}\label{DualQuadraticCone}
\Q^*(A)& =& \{V\in\S^n\ :\ U\cdot V\geq 0\,,\ \ U\in\Q(A)\} \\
\nonumber
& =& \{V\in\S^n\ :\ (gh\t+hg\t)\cdot V\ \geq\ 0
\,,\ \ g,h\in\G(A)\,,\ g\neq h\,,\ g\circ h\geq 0\} \\
\nonumber
& =& \{V\in\S^n\ :\ g\t Vh\ \geq\ 0
\,,\ \ g,h\in\G(A)\,,\ g\neq h\,,\ g\circ h\geq 0\}.
\end{eqnarray}
\ed

We are also interested is the following cone of diagonals
of matrices in $\Q(A)$.
\bd{DiagonalGraverCone}
The {\em diagonal Graver cone} of $A$ is the
cone of nonnegative vectors
$$\D(A)\ :=\ \cone\left\{g\circ h\ :\ g,h\in\G(A)\,,\
g\neq h\,,\ g\circ h\geq 0\right\}\ \subseteq\ \R_+^n\ .$$
The {\em dual diagonal Graver cone}
is its dual $\D^*(A)$ in $\R^n$ given by
\begin{eqnarray}\label{DualDiagonalCone}
\D^*(A)\ =\  \{v & : &  u\t v\geq 0\,,\ \ u\in\D(A)\} \\
\nonumber
\ =\ \{v & : & (g\circ h)\t v\ \geq\ 0
\,,\ \ g,h\in\G(A) \,,\ g\neq h\,,\ g\circ h\geq 0\} \\
\nonumber
\ =\ \{v & : & \sum g_ih_iv_i\ \geq\ 0
\,,\ \ g,h\in\G(A) \,,\ g\neq h\,,\ g\circ h\geq 0\}.
\end{eqnarray}
\ed

The following lemma provides some basic relations among the
above cones and more. All inclusions can be strict,
as is demonstrated in Examples \ref{LeftStrict} and
\ref{RightStrict} below.
In particular, it is interesting to note that $\D(A)$ is
the diagonal projection of $\Q(A)$, but $\D^*(A)$ is generally
strictly contained in the diagonal projection of $\Q^*(A)$.
\bl{DualCones}
The quadratic and diagonal Graver cones and their duals satisfy
\begin{eqnarray}\label{Relations}
\R_+^n\ \supseteq & \D(A) & = \
\left\{\diag(U)\ :\ U\in\Q(A)\right\}\ \supseteq\
\left\{u\ :\ \Diag(u)\in\Q(A)\right\}, \\
\nonumber
\R_+^n\ \subseteq & \D^*(A) & = \
\left\{v\ :\ \Diag(v)\in\Q^*(A)\right\}\ \subseteq\
\left\{\diag(V)\ :\ V\in\Q^*(A)\right\}.
\end{eqnarray}
\el
\bpr
First, $\D(A)\subseteq\R_+^n$ because it is generated by
nonnegative vectors. Therefore $\D^*(A)\supseteq(\R_+^n)^*=\R_+^n$.
To establish the top equality note that the following are equivalent:
$u\in\D(A)$; $u=\sum_k \mu_k (g_k\circ h_k)$ for some suitable
$\mu_k\geq 0$, $g_k,h_k\in\G(A)$; $u=\diag(U)$ with
$U={1\over 2}\sum_k\mu_k(g_k\otimes h_k+h_k\otimes g_k)$;
and $u=\diag(U)$ with $U\in\Q(A)$. To establish the bottom
equality note that the following are equivalent: $v\in\D^*(A)$;
$(g\circ h)\t v\geq 0$ for all suitable $g,h\in\G(A)$;
$V=\Diag(v)$ with $g\t V h\geq 0$ for all $g,h$;
and $V=\Diag(v)$ with $V\in\Q^*(A)$. The two remaining inclusions
on the right-hand sides follow from $\diag(\Diag(x))=x$.
This completes the proof of the lemma.
\epr
The next two examples show that all inclusions
in Lemma \ref{DualCones} can be strict.

\be{LeftStrict}
Consider the zero $1\times n$ matrix $A:=0$, whose Graver
basis is given by $\G(A)=\{\pm{\bf 1}_i\,:\,i=1,\dots,n\}$.
Then $g\circ h=0$ is the zero vector for all distinct
$g,h\in\G(A)$ in the same orthant. So the diagonal Graver cone
and its dual are $\D(A)=\{0\}\subsetneq\R_+^n$ and
$\D^*(A)=\R^n\supsetneq\R_+^n$ so the left inclusions
in (\ref{Relations}) are strict.
\ee
\be{RightStrict}
Consider the $1\times 3$ matrix $A:=(1\,\ 1\,\ 1)$ with
Graver basis $\G(A)=\pm\{(1,-1,0),(1,0,-1),(0,1,-1)\}$.
The quadratic Graver cone and its dual satisfy
\begin{eqnarray*}
\Q(A)\ =\ \cone\left\{
\left(\begin{array}{rrr}
  2 & -1 & -1 \\
 -1 &  0 &  1 \\
 -1 &  1 &  0 \\
\end{array}\right),
\left(\begin{array}{rrr}
  0 & -1 &  1 \\
 -1 &  2 & -1 \\
  1 & -1 &  0 \\
\end{array}\right),
\left(\begin{array}{rrr}
  0 &  1 & -1 \\
  1 &  0 & -1 \\
 -1 & -1 &  2 \\
\end{array}\right)
\right\}\ ,
\end{eqnarray*}
\begin{eqnarray}\label{IndefiniteExample}
\Q^*(A) & = &
\nonumber
\left\{\left(\begin{array}{rrr}
  a &  d & e \\
  d &  b & f \\
  e &  f & c \\
\end{array}\right)
\ :\
\begin{array}{rrr}
  a -d-e+f &  \geq & 0 \\
  b -d+e-f &  \geq & 0 \\
  c +d-e-f &  \geq & 0 \\
\end{array}\right\} \\
& \supseteq &
\left\{\left(\begin{array}{rrr}
  2a  &  a+b & a+c \\
  a+b &  2b  & b+c \\
  a+c &  b+c & 2c  \\
\end{array}\right)
\ :\ a,b,c\in\R\right\}.
\end{eqnarray}
The diagonal Graver cone and its dual are
$\D(A)=\R_+^n$ and $\D^*(A)=\R_+^n$.
Therefore, the top and bottom inclusions on the
right-hand side of equation (\ref{Relations}) are strict,
$$\D(A)=\R_+^n\ \supsetneq\ \{0\}\ =\
\left\{u\ :\ \Diag(u)\in\Q(A)\right\},$$
$$\D^*(A)=\R_+^n\ \subsetneq\ \R^n\ =\
\left\{\diag(V)\ :\ V\in\Q^*(A)\right\}.$$
\ee

\section{Quadratic integer minimization}

We proceed to establish our algorithmic
Theorems \ref{QuadraticIntegerMinimization}
and \ref{DiagonalIntegerMinimization}.
We focus on the situation of finite feasible sets,
which is natural in most applications.
But we do allow the lower and upper bounds $l,u\in\Zinf^n$ to have
infinite components for flexibility of modeling (for instance, it is
quite common in applications to have $l_i=0$ and $u_i=\infty$ for
all $i$, with the resulting feasible set typically still finite).
We also require our algorithms to identify and properly stop
when the set is infinite. So in all algorithmic statements,
an algorithm is said to {\em solve} a (nonlinear) discrete
optimization problem, if for every input,
it either finds an optimal solution, or asserts that
the problem is infeasible or the feasible set is infinite.
We begin with a simple lemma that shows that we can quickly
minimize a given quadratic function in a given direction.
\bl{GreedyAugmentation}
There is an algorithm that, given bounds $l,u\in\Z_{\infty}^n$,
direction $g\in\Z^n$, point $z\in\Z^n$ with $l\leq z\leq u$,
and quadratic function $f(x)=x\t Vx+w\t x+a$ with
$V\in\Z^{n\times n}$, $w\in\Z^n$, and $a\in\Z$,
solves in polynomial time the univariate problem
\begin{equation}\label{GreedyAugmentationEquation}
\min\{f(z+\mu g)\ :\ \mu\in\Z_+\,,\ l\leq z+\mu g\leq u\}.
\end{equation}
\el
\bpr
Let $S:=\{\mu\in\Z_+:l\leq z+\mu g\leq u\}$, and let $s:=\sup S$
which is easy to determine. If $s=\infty$ then
we conclude that $S$ is infinite and stop. Otherwise we need to
minimize the univariate quadratic function
$h(\mu):=f(z+\mu g)=h_2\mu^2+h_1\mu+h_0$ with $h_2:=g\t Vg$,
$h_1:=z\t Vg+g\t Vz+w\t g$, and $h_0:=z\t Vz+w\t z+a$ over
$S=\{0,1,\dots,s\}$. If $h_2\leq 0$, then $h$ is concave, and
the minimum over $S$ is attained at $\mu=0$ or $\mu=s$.
If $h_2>0$ then $h$ is convex with real minimum at
$\mu^*:=-{h_1\over 2h_2}$. Then minimizing $h$ over $S$
reduces to minimizing $h$ over
$S\cap\{0,\lfloor\mu^*\rfloor,\lceil \mu^*\rceil,s\}$.
\epr
A finite sum $u:=\sum_i v_i$ of vectors in $\R^n$
is called {\em conformal} if $v_i\sqsubseteq u$ for all $i$, and
hence all summands lie in the same orthant.
The following lemma shows that quadratic $f$ with
defining matrix in the dual quadratic Graver cone
is supermodular on conformal sums of nonnegative combinations
of elements of the Graver basis.
\bl{SeparableConvexConformal}
Let $A$ be any integer $m\times n$ matrix with quadratic Graver
cone $\Q(A)$. Let $f:\R^n\rightarrow\R$ be any quadratic
function $f(x)=x\t Vx+w\t x+a$ with $V\in\Q^*(A)$.
Let $x\in\R^n$ be any point, and let $\sum\mu_i g_i$ be any
conformal sum in $\R^n$ with $g_i\in\G(A)$ distinct elements in the
Graver basis of $A$ and $\mu_i\geq 0$ nonnegative scalars. Then
$$\Delta\ :=\ \left(f\left(x+\sum \mu_i g_i\right)-f(x)\right)
\ -\ \sum \left(f\left(x+\mu_i g_i\right)-f(x)\right)\ \geq\ 0.$$
\el
\bpr
We have
\begin{eqnarray*}
\!f\left(x+\sum \mu_i g_i\right)-f(x)=
\sum x\t V\mu_j g_j +\sum \mu_i g_i\t Vx+
\sum_{i,j} \mu_i g_i\t V \mu_j g_j + \sum w\t \mu_i g_i\ ,
\end{eqnarray*}
and
\begin{eqnarray*}
\sum \left(f(x+\mu_i g_i)-f(x)\right)& =&
\sum\left(x\t V\mu_i g_i+\mu_i g_i\t Vx+
\mu_i g_i\t V\mu_i g_i+w\t \mu_i g_i\right).
\end{eqnarray*}
Therefore we obtain
\begin{eqnarray*}
\Delta\ =\
\sum_{i,j} \mu_i g_i\t V \mu_j g_j-\sum \mu_i g_i\t V\mu_i g_i
\ =\ \sum_{i\neq j} \mu_i g_i\t V\mu_j g_j\ =\
\sum_{i\neq j}\mu_i\mu_j g_i\t V g_j\ \geq 0,
\end{eqnarray*}
because $g_i,g_j\in\G(A)$ satisfy $g_i\circ g_j\geq 0$
and $g_i\neq g_j$ for $i\neq j$, and $V$ is in $\Q^*(A)$.
\epr
We need two more useful properties of Graver bases. First we need
the following integer analogue of Carath\'eodory's theorem
of \cite{Seb} which we state without proof.
\bl{ConformalGraverSum}
Let $A$ be an integer $m\times n$ matrix, and let $\G(A)$
be its Graver basis. Then every $x\in\L^*(A)$ is a
conformal sum $x=\sum_{i=1}^t\mu_i g_i$ that involves
$t\leq 2n-2$ Graver basis elements $g_i\in\G(A)$ and
nonnegative integer coefficients $\mu_i\in\Z_+$.
\el
The next lemma provides a Graver basis criterion
for finiteness of integer programs.
\bl{GraverFiniteness}
Let $\G(A)$ be the Graver basis of matrix $A$, and let
$l,u\in\Z_{\infty}^n$. If there is some $g\in\G(A)$
satisfying $g_i\leq 0$ whenever $u_i<\infty$ and $g_i\geq 0$
whenever $l_i>-\infty$ then every set of the form
$S:=\{x\in\Z^n\,:\,Ax=b\,,\ l\leq x\leq u\}$ is either
empty or infinite, whereas if there is no such $g$, then every
set $S$ of this form is finite. Clearly, given the Graver basis,
the existence of such $g$ can be checked in polynomial time.
\el
\bpr
Suppose there is such $g$ and consider such $S$ containing
a point $x$. Then for all $\lambda\in\Z_+$ we
have $l\leq x+\lambda g\leq u$ and $A(x+\lambda g)=Ax=b$, and hence
$x+\lambda g\in S$ so $S$ is infinite. Next suppose $S$ is infinite.
Then $P:=\{x\in\R^n:Ax=b,\, l\leq x\leq u\}$
is unbounded, and hence has a recession vector, which we may assume
is integer, that is, a nonzero $h$ such that $x+\alpha h\in P$
for all $x\in P$ and $\alpha\geq 0$. Then $h\in\L^*(A)$ and
$h_i\leq 0$ whenever $u_i<\infty$ and $h_i\geq 0$ whenever
$l_i>-\infty$. By Lemma \ref{ConformalGraverSum}, the vector $h$
is a conformal sum $h=\sum g_i$ of vectors $g_i\in\G(A)$, each of
which also satisfies $g_i\leq 0$ whenever $u_i<\infty$ and
$g_i\geq 0$ whenever $l_i>-\infty$, providing such $g$.
\epr
Next we prove the main lemma underlying our algorithm, which shows
that, given the Graver basis, and an initial feasible point,
we can minimize a quadratic function with defining matrix in
the dual quadratic Graver cone in polynomial time.
\bl{QuadraticMinimizationLemma}
There is an algorithm that, given $A\in\Z^{m\times n}$,
its Graver basis $\G(A)$, bounds $l,u\in\Z_{\infty}^n$, point
$z\in\Z^n$ with $l\leq z\leq u$, and quadratic
$f(x)=x\t Vx+w\t x+a$ with integer $V\in\Q^*(A)$, $w\in\Z^n$,
and $a\in\Z$, solves in polynomial time the program
\begin{equation}\label{WeakSeparableConvexEquation}
\min\{f(x)=x\t Vx+w\t x+a\ :\
x\in\Z^n\,,\ Ax=b\,,\ l\leq x\leq u\}\ ,\quad b:=Az\ .
\end{equation}
\el
\bpr
First, apply the algorithm of Lemma \ref{GraverFiniteness}
to $\G(A)$ and $l,u$ and either detect that the feasible
set is infinite and stop, or conclude it is finite and continue.
Next produce a sequence of feasible points $x_0,x_1,\ldots,x_s$
with $x_0:=z$ the given input point, as follows.
Having obtained $x_k$, solve the minimization problem
\begin{equation}\label{GreedyEquation}
\min\{f(x_k+\mu g)\ :\ \mu\in\Z_+\,,\
g\in \G(A)\,,\ l\leq x_k+\mu g\leq u\,\}
\end{equation}
by applying the algorithm of Lemma \ref{GreedyAugmentation}
for each $g\in\G(A)$. If the minimal value in (\ref{GreedyEquation})
satisfies $f(x_k+\mu g)<f(x_k)$ then set $x_{k+1}:=x_k+\mu g$
and repeat, else stop and output the last point $x_s$ in the
sequence. Now, $Ax_{k+1}=A(x_k+\lambda g)=Ax_k=b$ by induction
on $k$, so each $x_k$ is feasible.
Because the feasible set is finite and the $x_k$ have decreasing
objective values and hence distinct, the algorithm terminates.

We now show that the point $x_s$ output by the algorithm is optimal.
Let $x^*$ be any optimal solution to (\ref{WeakSeparableConvexEquation}).
Consider any point $x_k$ in the sequence, and suppose that it is not optimal.
We claim that a new point $x_{k+1}$ will be produced and will satisfy
\begin{equation}\label{FastConvergence}
f(x_{k+1})-f\left(x^*\right)\ \leq\ {2n-3\over 2n-2}\left(f(x_k)-f(x^*)\right).
\end{equation}
By Lemma \ref{ConformalGraverSum}, we can write the difference
$x^*-x_k=\sum_{i=1}^t\mu_i g_i$ as conformal sum involving
$1\leq t\leq 2n-2$ elements $g_i\in \G(A)$ with all $\mu_i\in\Z_+$.
By Lemma \ref{SeparableConvexConformal},
$$f(x^*)-f\left(x_k\right)\ =\
f\left(x_k+\sum_{i=1}^t\mu_i g_i\right)-f(x_k)\ \geq\
\sum_{i=1}^t\left(f\left(x_k+\mu_i g_i\right)-f(x_k)\right).$$
Adding $t\left(f(x_k)-f(x^*)\right)$ on both
sides and rearranging terms, we obtain
$$\sum_{i=1}^t\left(f\left(x_k+\mu_i g_i\right)-f(x^*)\right)
\ \leq\ (t-1)\left(f(x_k)-f(x^*)\right).$$
Therefore there is some summand on the left-hand side satisfying
$$f\left(x_k+\mu_i g_i\right)-f(x^*)
\ \leq\ {t-1\over t}\left(f(x_k)-f(x^*)\right)
\ \leq\ {2n-3\over 2n-2}\left(f(x_k)-f(x^*)\right).$$
So the point $x_k+\mu g$ attaining
minimum in (\ref{GreedyEquation}) satisfies
$$f(x_k+\mu g)-f(x^*)\ \leq\ f\left(x_k+\mu_i g_i\right)-f(x^*)
\ \leq\ {2n-3\over 2n-2}\left(f(x_k)-f(x^*)\right),$$
and so indeed $x_{k+1}:=x_k+\mu g$ will be produced and will
satisfy (\ref{FastConvergence}). This shows that the last point
$x_s$ produced and output by the algorithm is indeed optimal.

We proceed to bound the number $s$ of points. Consider any $i<s$
and the intermediate non-optimal point $x_i$ in the sequence
produced by the algorithm. Then $f(x_i)>f(x^*)$ with both values
integer, and so repeated use of (\ref{FastConvergence}) gives
\begin{eqnarray*}
1\leq f(x_i)-f(x^*) & = & \prod_{k=0}^{i-1}
{{f(x_{k+1})-f(x^*)}\over{f(x_k)-f(x^*)}}\left(f(x)-f(x^*)\right) \\
& \leq & \left({2n-3\over 2n-2}\right)^i\left(f(x)-f(x^*)\right),
\end{eqnarray*}
and therefore
$$i\ \leq\ \left(\log{2n-2\over 2n-3}\right)^{-1}\log\left(f(x)-f(x^*)\right).$$
Therefore the number $s$ of points produced by the algorithm is at most one unit
larger than this bound, and using a simple bound on the logarithm, we obtain
$$s\ =\ O\left(n \log(f(x)-f(x^*))\right).$$
Thus, the number of points produced and the total running time are polynomial.
\epr
Next we show that, given the Graver basis, we can also find
an initial feasible point for assert that the given set is
empty or infinite, in polynomial time.
\bl{Feasibility}
There is an algorithm that, given integer $m\times n$ matrix $A$,
its Graver basis $\G(A)$, $l,u\in\Z_{\infty}^n$, and $b\in\Z^m$,
in polynomial time, either finds a feasible point in the set
$S\,:=\,\{x\in\Z^n\,:\,Ax=b,\, l\leq x\leq u\}$
or asserts that $S$ is empty or infinite.
\el
\bpr
Assume that $l\leq u$ and that $l_j<\infty$ and $u_j>-\infty$
for all $j$, because otherwise there is no feasible point.
Also assume that there is no $g\in\G(A)$ satisfying $g_j\leq 0$
whenever $u_j<\infty$ and $g_j\geq 0$ whenever $l_j>-\infty$, because
otherwise $S$ is empty or infinite by Lemma \ref{GraverFiniteness}.
Now, either detect there is no integer solution
to the system of equations $Ax=b$ (without the lower
and upper bound constraints) and stop, or determine some such
solution $x\in\Z^n$ and continue; it is well known that
this can be done in polynomial time, say, using the Hermite
normal form of $A$, see \cite{Sch}. Let
$$I\ :=\ \{j\,:\,l_j\leq x_j\leq u_j\}\ \subseteq\ \{1,\dots,n\}$$
be the set of indices of entries of $x$ that
satisfy their lower and upper bounds.
While $I\subsetneq \{1,\dots,n\}$ repeat the following procedure.
Pick any index $i\notin I$. Then either $x_i<l_i$ or $x_i>u_i$.
We describe the procedure only in the former case,
the latter being symmetric.
Update the lower and upper bounds by setting
$${\hat l}_j\ :=\ \min\{l_j,x_j\}\,,\quad\quad
{\hat u}_j\ :=\ \max\{u_j,x_j\}\,,\quad\quad j=1,\dots,n\ .$$
Solve in polynomial time the following linear integer program,
for which $x$ is feasible,
\begin{equation}\label{auxiliary_program}
\max\{z_i\ :\ z\in\Z^n\,,\ Az=b
\,,\ {\hat l}\leq z\leq {\hat u}\,,\ z_i\leq u_i\},
\end{equation}
by applying the algorithm of Lemma \ref{QuadraticMinimizationLemma}
using the function $f(z):=z\t 0z+{\bf 1}_i\t z+0$ with $V=0$
the zero matrix which is always in $\Q^*(A)$.
Now ${\hat l}_j>-\infty$ if and only if $l_j>-\infty$,
and ${\hat u}_j<\infty$ if and only if $u_j<\infty$.
So there is no $g\in\G(A)$ satisfying $g_j\leq 0$ whenever
${\hat u}_j<\infty$ and $g_j\geq 0$ whenever ${\hat l}_j>-\infty$,
and hence the feasible set of (\ref{auxiliary_program}) is finite
by Lemma \ref{GraverFiniteness} and has an optimal solution $z$.
If $z_i<l_i$ then assert that the set $S$ is empty and stop.
Otherwise, set $x:=z$, $I:=\{j\,:\,l_j\leq x_j\leq u_j\}$,
and repeat. Note that in each iteration, the cardinality of $I$
increases by at least one. Therefore, after at most $n$ iterations,
either the algorithm detects infeasibility, or $I=\{1,\dots,n\}$
is obtained, in which case the current point $x$ is feasible.
\epr
We are now in position to establish our theorem.

\vskip.2cm\noindent{\bf Theorem \ref{QuadraticIntegerMinimization}}
{\it There is an algorithm that, given $A\in\Z^{m\times n}$,
its Graver basis $\G(A)$, bounds $l,u\in\Z_{\infty}^n$,
$b\in\Z^m$, integer matrix $V\in\Q^*(A)$ in the dual quadratic
Graver cone, $w\in\Z^n$, and $a\in\Z$,
solves in polynomial time the quadratic integer program
\begin{equation*}
\min\{x\t Vx+w\t x+a\ :\ x\in\Z^n\,,\ Ax=b\,,\ l\leq x\leq u\}.
\end{equation*}
}
\bpr
Use the algorithm underlying Lemma \ref{Feasibility} to either
detect that the problem is infeasible or that the feasible set
is infinite and stop, or obtain a feasible point and use the
algorithm underlying Lemma \ref{QuadraticMinimizationLemma}
to obtain an optimal solution.
\epr
An important immediate consequence of Theorem
\ref{QuadraticIntegerMinimization} is that we can
efficiently minimize separable quadratic
functions defined by vectors in the dual diagonal Graver cone.
In particular, it applies to every convex separable
quadratic function (which can also be deduced from
the results of \cite{HOW} on separable convex functions).

\vskip.2cm\noindent{\bf Theorem \ref{DiagonalIntegerMinimization}}
{\it There is an algorithm that, given $A\in\Z^{m\times n}$,
its Graver basis $\G(A)$, bounds $l,u\in\Z_{\infty}^n$,
$b\in\Z^m$, integer vector $v\in\D^*(A)$ in the dual diagonal
Graver cone, and $w,a\in\Z^n$, solves in polynomial time the
separable quadratic program
\begin{equation*}
\min\{\sum_{i=1}^n\left(v_i x_i^2+w_ix_i+a_i\right)\ :\
x\in\Z^n\,,\ Ax=b\,,\ l\leq x\leq u\}.
\end{equation*}
In particular, this applies to any convex separable
quadratic, that is, with $v\in\Z_+^n$.
}
\vskip.2cm
\bpr
First, for any $v\in\D^*(A)$ we have $V:=\Diag(v)\in\Q^*(A)$
by Lemma \ref{DualCones}.
Hence, by Theorem \ref{QuadraticIntegerMinimization},
we can minimize in polynomial time the quadratic function
$$\sum_{i=1}^n\left(v_i x_i^2+w_ix_i+a_i\right)
\ =\ x\t Vx+w\t x+\sum_{i=1}^n a_i\ .$$

Second, if the separable quadratic function is convex, which
is equivalent to its defining vector $v$ being nonnegative,
then $v\in\R_+^n\subseteq\D^*(A)$ by Lemma \ref{DualCones}
again. Hence the second statement of the theorem now
follows from the first statement.
\epr

\section{Nonconvex solvable quadratics and matroids}

Consider the quadratic minimization problem,
with the Graver basis of $A$ given,
\begin{equation}\label{PSD}
\min\left\{f(x)=x\t Vx+w\t x+a\ :\ x\in\Z^n\,,
\ Ax=b\,,\ l\leq x\leq u\right\}.
\end{equation}
The function $f$ is convex if and only if its defining
matrix $V$ is positive semidefinite, that is, if $x\t Vx\geq 0$
for all $x\in\R^n$. Let $\S_+^n\subset\S^n$ denote the cone
of symmetric positive semidefinite matrices.
Now, on the one hand, if $V\in\Q^*(A)$ then, by
Theorem \ref{QuadraticIntegerMinimization}, we can solve
problem (\ref{PSD}) efficiently. On the other hand, if $V\in\S_+^n$
then $f$ is convex, and problem (\ref{PSD}) may seem to be easier,
but remains NP-hard even for rank-$1$ matrices $V=vv\t\in\S_+^n$
by Proposition \ref{QuadraticHardness}.
So it is unlikely that $\Q^*(A)$ contains $\S_+^n$, and it is
interesting to consider the relation between these matrix cones.

For this, we need a couple of basic facts about positive
semidefinite matrices. First, Note that for any vector
$u\in\R^n$, the rank-$1$ matrix $uu\t$ is in $\S_+^n$ because
$x\t(uu\t)x=(u\t x)^2\geq 0$ for all $x\in\R^n$, whereas
for any two linearly independent vectors $g,h\in\R^n$,
the rank-$2$ matrix $gh\t+hg\t$ is in $\S^n\setminus\S_+^n$
because there is an $x\in\R^n$ with $g\t x=1$ and $h\t x=-1$
and hence $x\t(gh\t+hg\t)x=2(g\t x)(h\t x)=-2<0$.
Second, the cone of symmetric positive semidefinite matrices
is self dual, that is, $(\S_+^n)^*=\S_+^n$. To see this,
note that if $U\in\S^n\setminus\S_+^n$ then there is an
$x\in\R^n$ with $(x\otimes x)\cdot U=x\t Ux<0$ so
$U\notin(\S_+^n)^*$; and if $V\in\S_+^n$ has rank $r$,
then $V=\sum_{i=1}^r x_i\otimes x_i$ for some $x_i\in\R^n$
and hence $U\cdot V=\sum_{i=1}^r x_i\t U x_i\geq 0$
for all $U\in\S_+^n$, so $V\in(\S_+^n)^*$.

So we can conclude the following. In the rare situation
where each orthant of $\R^n$ contains at most one element of
$\G(A)$, we have $\Q(A)=\{0\}$
and $\Q^*(A)=\S^n$, so Theorem \ref{QuadraticIntegerMinimization}
enables to solve problem (\ref{PSD}) for any quadratic function.
In the more typical situation, where some orthant does contain two
elements $g,h\in\G(A)$, the corresponding generator of $\Q(A)$
satisfies $gh\t+hg\t\in\S^n\setminus\S_+^n$ and hence
$\Q(A)\nsubseteq\S_+^n$. By self duality of $\S_+^n$,
we obtain $\S_+^n=(\S_+^n)^*\nsubseteq\Q^*(A)$.
So we cannot solve problem (\ref{PSD}) for all convex quadratics,
reflecting the NP-hardness of the convex problem.
But we do typically also have $\Q^*(A)\nsubseteq\S_+^n$,
that is, we can solve problem (\ref{PSD}) in polynomial time
for various nonconvex quadratics. For instance, in
Example \ref{RightStrict}, the matrix in (\ref{IndefiniteExample})
is not positive semidefinite for all $a,b,c<0$. Moreover,
by Lemma \ref{DualCones},
$\R_+^n\subseteq\D^*(A)=\{v\,:\,\Diag(v)\in\Q^*(A)\}$,
so $\Q^*(A)\setminus\S_+^n\neq\emptyset$ whenever
$\D^*(A)\setminus\R_+^n\neq\emptyset$.

We proceed to discuss this diagonal case, where the
function $f$ is defined by a diagonal matrix $V=\Diag(v)$
for some $v\in\R^n$, that is, $f$ is separable of the form
$f(x)=\sum_i (v_ix_i^2+w_ix_i+a_i)$. In this case, $f$ is
convex if and only if $v$ is nonnegative. As noted in
Lemma \ref{DualCones}, the dual diagonal Graver cone $\D^*(A)$
always contains the nonnegative orthant $\R_+^n$. We proceed
to characterize those matrices $A$ for which this inclusion
is strict, so that $\D^*(A)\setminus\R_+^n\neq\emptyset$
and Theorem \ref{DiagonalIntegerMinimization}
enables to solve problem (\ref{PSD}) in polynomial time also for
various nonconvex separable quadratics.

For this we need a few more
definitions. A {\em circuit} of an integer matrix $A$ is an
element $c\in\L^*(A)$ whose support $\supp(c)$ is minimal under
inclusion and whose entries are relatively prime. We denote the
set of circuits of $A$ by $\C(A)$. It is easy to see that for every
integer matrix $A$, the set of circuits is contained in the Graver
basis, that is, $\C(A)\subseteq\G(A)$. Recall that a finite sum
$u:=\sum_i v_i$ of vectors in $\R^n$ is {\em conformal} if
$v_i\sqsubseteq u$ for all $i$, and hence all summands lie in the
same orthant. The following property of circuits is well known.
For a proof see, for instance, \cite{Onn} or \cite{Stu}.
\bl{ConformalCircuitSum}
Let $A$ be an integer matrix. Then every $x\in\L^*(A)$ is a
conformal sum $x=\sum_i\alpha_i c_i$ involving circuits
$c_i\in\C(A)$ and nonnegative real coefficients $\alpha_i\in\R_+$.
\el
It turns out that the matroid of linear dependencies
on the columns of the integer $m\times n$ matrix $A$
(over the reals or integers) plays a central role in the
characterization we are heading for. A {\em matroid-circuit}
is any set $C\subseteq\{1,\dots,n\}$ that is the support
$C=\supp(c)$ of some circuit $c\in\C(A)$ of $A$.
Note that a circuit $c$ is in $\C(A)$ if and only if its
antipodal $-c$ is, and if $c,e\in\C(A)$ are circuits with
$c\neq\pm e$ then $\supp(c)\neq \supp(e)$. We denote the set
of matroid-circuits of $A$, that is, the set of supports of
circuits in $\C(A)$, by $\M(A):=\{\supp(c)\,:\,c\in\C(A)\}$,
and refer to it simply as the {\em matroid} of $A$.
For instance, for the $1\times 3$ matrix $A:=(1\,\ 2\,\ 1)$ we have
$$\C(A)\ =\ \pm\left\{(2,-1,0),(0,-1,2),(1,0,-1)\right\},\quad
\M(A)\ =\ \left\{\{1,2\},\{2,3\}),\{1,3\}\right\}.$$

We now characterize those matrices $A$ for which
$\D^*(A)$ strictly contains $\R_+^n$.
\bt{Characterization}
The dual diagonal Graver cone of every
integer $m\times n$ matrix $A$ satisfies $\D^*(A)\supseteq\R_+^n$,
and the inclusion is strict if and only if there is
$1\leq k\leq n$ such that $C\cap E\neq\{k\}$ for every
two distinct matroid-circuits $C,E\in\M(A)$ of $A$.
\et
\bpr
We prove the dual statements about the diagonal Graver cone.
By definition $\D(A)\subseteq\R_+^n$, and the inclusion is strict
if and only if some unit vector ${\bf 1}_k$ is not in $\D(A)$.
Therefore it suffices to prove that, for any $1\leq k\leq n$,
we have ${\bf 1}_k\in\D(A)$ if and only if there are two
distinct matroid-circuits $C,E\in\M(A)$ with $C\cap E=\{k\}$.

Suppose first $C,E\in\M(A)$ are distinct matroid-circuits with
$C\cap E=\{k\}$. Then there are $c,e\in\C(A)$ with $c\neq\pm e$
such that $\supp(c)=C$ and $\supp(e)=E$. Replacing $e$ by
$-e\in\C(A)$ if necessary we may assume that $c_ke_k>0$.
Then $c\circ e\geq 0$, $c\neq e$, and $c,e\in\G(A)$ imply that
$c\circ e=c_ke_k{\bf 1}_k$ is a generator of $\D(A)$,
and hence ${\bf 1}_k\in\D(A)$. Conversely, suppose
${\bf 1}_k\in\D(A)$. Because $\D(A)\subseteq\R_+^n$, some
nonnegative multiple of ${\bf 1}_k$ must be one of the generators.
So there are $g,h\in\G(A)$ with $g\circ h\geq 0$ and $g\neq h$
such that $g\circ h$ is a nonnegative multiple of ${\bf 1}_k$,
and hence $\supp(g)\cap\supp(h)=\{k\}$.
By Lemma \ref{ConformalCircuitSum} we have $g=\sum_i\alpha_i c_i$
and $h=\sum_j\alpha_j e_j$ conformal sums of circuits
with nonnegative coefficients. Then $\supp(g)=\cup\,\supp(c_i)$ and
$\supp(h)=\cup\,\supp(e_j)$, and hence there are $c_i$ and $e_j$
among these circuits such that $\supp(c_i)\cap\supp(e_j)=\{k\}$.
Let $C:=\supp(c_i)$ and $E:=\supp(e_j)$ be the corresponding
matroid-circuits of $A$. It remains to show
that $C$ and $E$ are distinct. Suppose indirectly that $C=E$.
Then $C=E=C\cap E=\{k\}$. This implies that the $k$-th column of $A$
is $0$ and $c_i=e_j=\pm{\bf 1}_k$. But then $c_i\sqsubseteq g$
and $e_j\sqsubseteq h$, and therefore $g=c_i=e_j=h$ which is a
contradiction. So $C\neq E$, and the proof is complete.
\epr
It is interesting to emphasize that the characterization in
Theorem \ref{Characterization} is in terms of only the matroid of
$A$, that is, the linear dependency structure on the columns of $A$.
The algorithm of Theorem \ref{DiagonalIntegerMinimization}
enables to solve in polynomial time the program
\begin{equation*}
\min\{\sum_{i=1}^n\left(v_i x_i^2+w_ix_i+a_i\right)\ :\
x\in\Z^n\,,\ Ax=b\,,\ l\leq x\leq u\}
\end{equation*}
for all separable quadratics with $v\in\D^*(A)$ and in particular
for all separable convex quadratic functions with $v\in\R_+^n$.
So the algorithm can solve the program moreover for some separable
nonconvex quadratic functions precisely when the matroid of $A$
satisfies the criterion of Theorem \ref{Characterization}.
Here are some concrete simple examples.
\be{SeparableBox}
Consider again Example \ref{LeftStrict} with $A:=0$ the
zero $1\times n$ matrix having
Graver basis $\G(A)=\{\pm{\bf 1}_i\,:\,i=1,\dots,n\}$.
Then the set of matroid-circuits of $A$ is
$\M(A)=\{\{1\},\dots,\{n\}\}$. Therefore $C\cap E=\emptyset$ for
all distinct $C,E\in\M(A)$ and the condition of
Theorem \ref{Characterization}
trivially holds, so $\D^*(A)\supsetneq\R_+^n$.
In fact, here $\D^*(A)=\R^n$.
\ee
\be{Digraphs}{\bf Directed graphs.}
Let $G$ be a directed graph, and let $A$ be its $V\times E$
incidence matrix, with $A_{v,e}:=1$ if vertex $v$ is the
head of directed edge $e$, $A_{v,e}:=-1$ if $v$ is the
tail of $e$, and $A_{v,e}:=0$ otherwise. The set $\M(A)$
of matroid-circuits consists precisely of all subsets
$C\subseteq E$ that are circuits of the undirected graph
underlying $G$. The set $\C(A)$ of circuits consists
of all vectors $c\in\{-1,0,1\}^E$ obtained from some
matroid circuit $C\subseteq E$ by choosing any of
its two orientations and setting $c_e:=1$ if directed edge
$e\in C$ agrees with the orientation, $c_e:=-1$ if $e$
disagrees, and $c_e:=0$ if $e\notin C$. The Graver basis
is equal to the set of circuits, $\G(A)=\C(A)$.
By Theorem \ref{Characterization} we have
$\D^*(A)\supsetneq\R_+^E$ if and only if there is an edge
$e\in E$ such that no two distinct circuits $C,C'$ of
the underlying undirected graph satisfy $C\cap C'=\{e\}$.
\ee
\be{MomentCurve}{\bf Generic Matrices.}
Let $A$ be a generic integer $m\times n$ matrix, that is,
a matrix for which every set of  $m$ columns is linearly independent,
say, the matrix defined by $A_{i,j}:=j^i$ for all $i,j$,
whose columns are distinct points on the moment curve in $\R^m$.
Then the matroid of $A$ is uniform, that is, its matroid-circuits
are exactly all $(m+1)$-subsets of $\{1,\dots,n\}$.
Suppose $n\leq 2m$.
Then every distinct $C,E\in\M(A)$ satisfy $|C\cap E|\geq 2$, and
hence $\D^*(A)\supsetneq\R_+^n$ by Theorem \ref{Characterization}.
So, by Theorem \ref{QuadraticIntegerMinimization},
\begin{equation*}
\min\{\sum_{i=1}^n\left(v_i x_i^2+w_ix_i+a_i\right)\ :\
x\in\Z^n\,,\ Ax=b\,,\ l\leq x\leq u\}
\end{equation*}
can be solved in polynomial time for all such $A$, all $b\in\Z^m$
and $l,u\in\Zinf^n$, all convex and some nonconvex separable
quadratic functions defined by data $v,w,a\in\Z^n$.
\ee

\section{Higher degree polynomial functions}

The algorithm that underlies our algorithmic
Theorem \ref{QuadraticIntegerMinimization}
using the Graver basis is conceptually quite simple.
First, it finds in polynomial time a feasible point.
Then it keeps improving points iteratively, as long as possible,
where, at each iteration, it takes the best possible improving
step attainable along any Graver basis element. It outputs the
last point from which no further Graver improvement is possible.

We now proceed to show that the results of the previous sections
can be extended to multivariate polynomials of higher,
arbitrary, degree. We will define a hierarchy of cones,
and whenever a polynomial function will lie in the corresponding
cone, the algorithm outlined above will converge to the optimal
solution in polynomial time.

It will be convenient now to make more extensive use of tensor
notation, and to work with the tensored, nonsymmetrized
form of a polynomial function. We use
$$\otimes_d\R^n\ := \R^n\otimes\cdots\otimes\R^n\ ,\quad
\otimes_d x\ :=\ x\otimes\cdots\otimes x\,,\ \ x\in\R^n$$
for the $d$-fold tensor product of $\R^n$ with itself
and for the rank-$1$ tensor that is the $d$-fold product of a
vector $x$ with itself, respectively. Note that the
$(i_1,\dots,i_d)$-th entry of $\otimes_d x$ is the product
$x_{i_1}\cdots x_{i_d}$ of the corresponding entries of $x$.
We denote the standard inner product on the tensor space by
$$\l U,V\r \ :=\ \sum_{i_1=1}^n\cdots \sum_{i_d=1}^n
U_{{i_1},\dots,{i_d}}V_{{i_1},\dots,{i_d}}
\ , \quad\quad U,V\ \in\ \otimes_d\R^n\ .$$
In particular, in the vector space $\R^n$ we have
$\l x,y\r =x\t y$ and in the matrix space
$\R^n\otimes \R^n$ we have $\l U,V\r =U\cdot V$.
Note that for any two rank-$1$ tensors we have
$$\l x^1\otimes \cdots \otimes x^d,
y^1\otimes \cdots \otimes y^d \r
\ =\ \prod_{k=1}^d \l x^k,y^k \r\ .$$

For simplicity, we restrict attention to homogeneous polynomials,
also termed {\em forms}. A form $f(x)$ of degree $d$ in the vector
of $n$ variables $x=(x_1,\dots,x_n)$ can be compactly defined by
a single tensor $F\in\otimes_d\R^n$ that collects
all coefficients, by
$$f(x)\ :=\ \l F,\otimes_d x\r\ =\
\sum_{i_1=1}^n\cdots \sum_{i_d=1}^n
F_{{i_1},\dots,{i_d}}x_{i_1}\cdots x_{i_d}\ .$$
For instance, the form $f(x)=(x_1+x_2+x_3)^3$
of degree $d=3$ in $n=3$ variables can be written as
$f(x)=\l F,\otimes_3 x\r=
\l\otimes_3{\bf 1},\otimes_3 x\r=\l{\bf 1},x\r^3$
with ${\bf 1}$ the all-ones vector in $\R^3$ and
$F=\otimes_3{\bf 1}$ the all-ones tensor in $\otimes_3\R^3$,
with $F_{i_1,i_2,i_3}=1$ for $i_1,i_2,i_3=1,2,3$.

Let $A$ be any integer $m\times n$ matrix, and let $\G(A)$
be its Graver basis. For each degree $d\geq 2$ we now define a
cone $\P_d(A)$ in the tensor space $\otimes_d\R^n$ as follows.
\bd{TensorGraverCone}
The {\em Graver cone of degree $d$} of an integer $m\times n$
matrix $A$ is the cone $\P_d(A)\subseteq\otimes_d\R^n$
generated by the rank-$1$ tensors $g^1\otimes \cdots \otimes g^d$
where the $g^i$ are elements of $\G(A)$ that lie in the same
orthant and are not all the same, that is
\begin{eqnarray*}
\hskip-.1cm
\P_d(A):=\cone\{g^1\otimes\cdots\otimes g^d\,:\,
g^i\in\G(A),\, g^i\circ g^j\geq 0\
\mbox{for all $i,j$},\,g^i\neq g^j\ \mbox{for some $i,j$}\}.
\end{eqnarray*}
The {\em dual Graver cone of degree $d$}
is its dual $\P_d^*(A)$ in $\otimes_d\R^n$ given by
\begin{eqnarray*}\label{DualTensorCone}
\P_d^*(A) = \{V\in\otimes_d\R^n \!\! & : \!\! &
\l U,V\r\geq 0,\ U\in\P_d(A)\}\ =\
\{V\, :\, \l g^1\otimes\cdots\otimes g^d,V\r\geq 0, \\
\nonumber
& &\ \  g^i\in\G(A),\ g^i\circ g^j\geq 0\ \
\mbox{for all $i,j$},\ g^i\neq g^j\ \mbox{for some $i,j$}\}.
\end{eqnarray*}
\ed
Note that $\P_2(A)$ is the nonsymmetrized version of $\Q(A)$,
that is, $\Q(A)=\P_2(A)\cap\S^n$.

One of the key ingredient in extending our algorithmic results
to polynomials of arbitrary degree is the following analogue of
Lemma \ref{SeparableConvexConformal} which establishes the
supermodularity of polynomial functions that lie in suitable cones.
We need one more piece of terminology.
Let $D:=\{1,\dots,d\}$ and for $0\leq k\leq d$ let
$D\choose k$ be the set of all $k$-subsets of $D$.
A $k$-dimensional {\em subtensor} of a $d$-dimensional tensor
$$F\ =\ (F_{i_1,\dots,i_d}\ :\
1\leq i_1,\dots,i_d\leq n) \ \in\ \otimes_d\R^n$$
is any of the ${d\choose k}n^{d-k}$ tensors
$T\in\otimes_k\R^n$
obtained from $F$ by choosing $I\in{D\choose k}$,
letting each index $i_j$ with $j\in I$ vary from $1$ to $n$,
and fixing each index $i_j$ with $j\notin I$ at some value
between $1$ and $n$. For instance, the $k$-dimensional tensor
obtained by choosing $I=\{1,\dots,k\}$ and fixing some values
$1\leq i_{k+1},\dots,i_d\leq n$ is
$$T\ =\ (T_{i_1,\dots,i_k}:=F_{i_1,\dots,i_k,i_{k+1},\dots,i_d}
\ :\ 1\leq i_1,\dots,i_k\leq n) \ \in\ \otimes_k\R^n\ .$$

For an integer $m\times n$ matrix $A$, let
$\K_d(A)\subseteq\otimes_d\R^n$ be the cone of those
tensors $F$ such that, for all $2\leq k\leq d$, every
$k$-dimensional subtensor of $F$ is in $\P_k^*(A)$.
\bl{Supermodularity}
Let $A$ be integer $m\times n$ matrix. Let
$f:\R^n\rightarrow\R$ be degree $d$ form given by
$f(x)=\l F,\otimes_d x\r$ with $F\in\K_d(A)$. Let $x\in\R_+^n$ be nonnegative
and $\sum_{r=1}^t\mu_r g^r$ conformal sum in $\R^n$ with
$g^r\in\G(A)$ distinct and $\mu_r\geq 0$ nonnegative scalars. Then
$$\Delta\ :=\ \left(f\left(x+\sum_{r=1}^t\mu_r g^r\right)-f(x)\right)\ -\
\sum_{r=1}^t \left(f\left(x+\mu_r g^r\right)-f(x)\right)\ \geq\ 0.$$
\el
\bpr
To simplify the derivation we assume that all $\mu_r=1$.
The same argument goes through in exactly the same way
for arbitrary nonnegative $\mu_r$. For $r=1,\dots,t$,
\begin{eqnarray*}
 f(x+g^r)-f(x)
& = & \l F,\otimes_d (x+g^r)\r - \l F,\otimes_d x\r \\
& = & \l F,g^r\otimes x\otimes\cdots\otimes x\r\ +\
\cdots\ +\ \l F,x\otimes\cdots\otimes x\otimes g^r\r \\
& + & \sum_{k=2}^d\sum\left\{\l F,u^1\otimes\cdots\otimes u^d\r\ :\
I\in{D\choose k}\,,\
u^i=\left\{
\begin{array}{ll}
  g^r, & i\in I \\
  x, & i\notin I
  \end{array}
  \right.
\right\}.
\end{eqnarray*}
Similarly,
\begin{eqnarray*}
f(x+\sum_{r=1}^t g^r)\!\!\!\! & -\!\!\!\! & f(x)
\ =\ \left\l F,\otimes_d \left(x+\sum_{r=1}^tg^r\right)\right\r
 - \l F,\otimes_d x\r \\
& = & \l F,\sum_{r=1}^tg^r\otimes x\otimes\cdots\otimes x\r \ +\
\cdots\ +\ \l F,x\otimes\cdots\otimes x\otimes \sum_{r=1}^tg^r\r \\
& + & \sum_{k=2}^d\sum\left\{\l F,u^1\otimes\cdots\otimes u^d\r\ :\
I\in{D\choose k}\,,\
u^i=\left\{
\begin{array}{ll}
  \sum_{r=1}^tg^r, & i\in I \\
  x, & i\notin I
  \end{array}
  \right.
\right\}.
\end{eqnarray*}
Therefore,
\begin{eqnarray}\label{Delta}
\Delta & =&
\sum_{k=2}^d\sum
\left\{\l F,u^1\otimes\cdots\otimes u^d\r
\ -\ \sum_{r=1}^t\l F,v^{r,1}\otimes\cdots\otimes v^{r,d}\r
\ :\ I\in{D\choose k}\,,\right.\\
\nonumber
& & \left. \hskip5cm u^i=\left\{
\begin{array}{ll}
  \sum_{r=1}^tg^r, & i\in I \\
  x, & i\notin I
  \end{array}
  \right.\,,\
v^{r,i}=\left\{
\begin{array}{ll}
  g^r, & i\in I \\
  x, & i\notin I
  \end{array}
  \right.
\right\}.
\end{eqnarray}
Now, consider any $2\leq k\leq d$ and any $I\in{D\choose k}$.
For simplicity of the indexation, we assume that $I=\{1,\dots,k\}$.
The derivation for other $I$ is completely analogous. For each
choice of indices $1\leq i_{k+1},\dots,i_d\leq n$ let
$T(i_{k+1},\dots,i_d)$ be the $k$-dimensional subtensor of $F$
obtained by letting $i_1,\dots,i_k$ vary and fixing
$i_{k+1},\dots,i_d$ as chosen. Then the corresponding summand
of $\Delta$ in the expression (\ref{Delta}) above satisfies
\begin{eqnarray}\label{Subtensor}
&&\hskip-1.5cm\left\l F,\otimes_k\left(\sum_{r=1}^tg^r\right)
\otimes(\otimes_{d-k}x)\right\r
\ -\ \sum_{r=1}^t\l F,(\otimes_k g^r)\otimes(\otimes_{d-k}x)\r \\
\nonumber
&=&\sum_{i_{k+1}=1}^n\cdots \sum_{i_d=1}^n
x_{i_{k+1}}\cdots x_{i_d}\left\l T(i_{k+1},\dots,i_d),
\otimes_k\left(\sum_{r=1}^tg^r\right)-
\sum_{r=1}^t\otimes_k g^r\right\r\ .
\end{eqnarray}
The summand in (\ref{Subtensor}) above which corresponds
to $1\leq i_{k+1},\dots,i_d\leq n$ satisfies
\begin{eqnarray}\label{Duality}
&&\hskip-.5cm\left\l T(i_{k+1},\dots,i_d),
\otimes_k\left(\sum_{r=1}^tg^r\right)-
\sum_{r=1}^t\otimes_k g^r\right\r\ =\ \\
\nonumber
&&\sum\left\{\left\l T(i_{k+1},\dots,i_d),g^{r_1}\otimes
\cdots\otimes g^{r_k}\right\r\ :\ 1\leq r_1,\dots,r_k\leq t
,\ \ \mbox{not all $r_i$ the same}\right\}.
\end{eqnarray}
Now, because all the $g^r$ are in the same orthant, and all
$k$-dimensional subtensors of $F$ lie in the dual Graver cone
$\P^*_k(A)$, each summand on the right-hand side of
(\ref{Duality}) above satisfies $\l T(i_{k+1},\dots,i_d),g^{r_1}
\otimes\cdots\otimes g^{r_k}\r\geq 0$, and so the left-hand side
of (\ref{Duality}) is nonnegative as well. Because $x\in\R_+^n$
is nonnegative, each summand on the right-hand side of
(\ref{Subtensor}) above is nonnegative, and so the left-hand
side of (\ref{Subtensor}) is nonnegative as well. Because this
holds for all $2\leq k\leq d$ and all $I\in{D\choose k}$,
we obtain that each summand on the right-hand side of (\ref{Delta})
is nonnegative, and so $\Delta\geq 0$ as claimed.
\epr
A second key ingredient is the following analogue of
Lemma \ref{GreedyAugmentation} which shows that we can efficiently
minimize a given form of any fixed degree $d$ in a given direction.
\bl{FormAugmentation}
For every fixed $d$, there is an algorithm that,
given $l,u\in\Z_{\infty}^n$, $z,g\in\Z^n$ with $l\leq z\leq u$,
and $f(x)=\l F,\otimes_d x\r$ with $F\in\otimes_d\Z^n$,
solves in polynomial time
\begin{equation}\label{GreedyAugmentationEquation}
\min\{f(z+\mu g)\ :\ \mu\in\Z_+\,,\ l\leq z+\mu g\leq u\}.
\end{equation}
\el
\bpr
Let $S:=\{\mu\in\Z_+:l\leq z+\mu g\leq u\}$, and let $s:=\sup S$
which is easy to determine. If $s=\infty$ then
we conclude that $S$ is infinite and stop. Otherwise we need to
minimize the univariate degree $d$ polynomial
$h(\mu):=\l F,\otimes_d(z+\mu g)\r=\sum_{i=0}^d h_i\mu^i$,
whose coefficients $h_i$ can be easily computed from $F$,
over $S=\{0,1,\dots,s\}$.

Outline: use repeated bisections
and Sturm's theorem which allows us to count the number of real
roots of $h$ in any interval using the Euclidean algorithm on
$h(\mu)=\sum_{i=0}^d h_i\mu^i$ and its derivative
$h'(\mu)=\sum_{i=0}^{d-1} (i+1)h_{i+1}\mu^i$,
to find intervals $[r_i,s_i]$, $i=1,\dots,d$ (possibly
with repetitions if $h$ has multiple roots) containing
each real root of $h$, and such that $s_i-r_i<1$ for all $i$.
Then minimizing $h$ over $S$ reduces to minimizing $h$ over
$S\cap\{0,\lceil r_1\rceil,\lfloor s_1\rfloor,
\dots \lceil r_d\rceil,\lfloor s_d \rfloor, s\}$.
\epr

We can now establish our theorem on polynomial integer minimization.

\vskip.2cm\noindent{\bf Theorem \ref{PolynomialIntegerMinimization}}
{\it For every fixed $d$ there is an algorithm that,
given integer $m\times n$ matrix $A$, its Graver basis $\G(A)$,
$b\in\Z^m$, and degree $d$ integer homogenous polynomial
$f(x)=\l F,\otimes_d x\r$ with $F\in\K_d(A)$, solves in
polynomial time the polynomial program
\begin{equation*}
\min\{f(x)=\l F,\otimes_dx\r\ :\ x\in\Z^n\,,\ Ax=b\,,\ x\geq 0\}.
\end{equation*}
}
\bpr
First, use the algorithm of Lemma \ref{Feasibility}
to either detect that the problem is infeasible or that the
feasible set is infinite and stop, or obtain a feasible point
and continue. Now, apply the algorithm of
Lemma \ref{QuadraticMinimizationLemma} precisely as it is,
using the given form $f(x)$ instead of a quadratic.
Lemmas \ref{Supermodularity} and \ref{FormAugmentation}
now assure that the analysis of this algorithm in the proof of
Lemma \ref{QuadraticMinimizationLemma} carries through
precisely as before, and guarantee that the algorithm will
find an optimal solution in polynomial time.
\epr

\vskip.6cm\noindent {\small Jon Lee}\newline
\emph{IBM T.J. Watson Research Center, Yorktown Heights, USA}\newline
\texttt{jonlee@us.ibm.com}

\vskip.3cm\noindent {\small Shmuel Onn}\newline
\emph{Technion - Israel Institute of Technology, Haifa, Israel}\newline
\texttt{onn@ie.technion.ac.il}

\vskip.3cm\noindent {\small Lyubov Romanchuk}\newline
\emph{Technion - Israel Institute of Technology, Haifa, Israel}\newline
\texttt{lyuba@techunix.technion.ac.il}

\vskip.3cm\noindent {\small Robert Weismantel} \newline
\emph{ETH, Z\"urich, Switzerland}\newline
\texttt{robert.weismantel@ifor.math.ethz.ch}

\end{document}